\documentclass[12pt]{amsart} 
\usepackage{amsmath, amssymb}
\usepackage[colorlinks=true]{hyperref}

\newcommand{\bbK}{{\mathbb K}}
\newcommand{\bbQ}{\mathbb{Q}}
\newcommand{\bbZ}{\mathbb{Z}}

\newcommand{\cC}{{\mathcal{C}}}
\newcommand{\cD}{{\mathcal{D}}}
\newcommand{\cN}{{\mathcal{N}}}
\newcommand{\cO}{{\mathcal{O}}}

\DeclareMathOperator{\core}{sf}
\DeclareMathOperator{\Gal}{Gal}

\newtheorem{lemma}{Lemma}

\newtheorem{theorem}{Theorem}

\theoremstyle{remark}
\newtheorem{note}{Note}

\begin{document}

\title[Thue's Fundamentaltheorem, II]{Thue's Fundamentaltheorem, II: Further Refinements and Examples}

\author{Paul M. Voutier}
\address{London, UK \\Paul.Voutier@gmail.com}

\begin{abstract}
In this paper, we sharpen and simplify our earlier results based on Thue's
Fundamentaltheorem and use it to obtain effective irrationality measures for
certain roots of particular polynomials of the form $(x-\sqrt{t})^{n}+(x+\sqrt{t})^{n}$,
where $n \geq 4$ is a positive integer and $t$ is a negative integer. For $n=4$
and $n=5$, we find infinitely many such algebraic numbers.
\end{abstract}


\maketitle

\section{Introduction}

In earlier papers \cite{Chen1, CV, LPV, TVW}, several authors have used Thue's
Fundamentaltheorem to completely solve several families of Thue equations and inequalities.
In \cite{Vout2, Vout3}, we simplified the statement of Thue's
Fundamentaltheorem and investigated the conditions under which it yields
effective irrationality measures for algebraic numbers.

In those papers, we attempted to simplify our statements by restricting
$d$ defined there to be a rational integer. However, this results in the need
for the quantities $g_{4}$ and $g_{5}$ in the definition of $g$ when the base
field is $\bbQ$ (see Corollary~3.7 of \cite{Vout3}). Furthermore, the results
are sometimes weaker than they need to be. By allowing $d$ to be the square root
of rational integer, we can both simplify and strengthen our previous results.

We use this new result to consider new examples as well. In particular,
roots of the polynomial
\begin{displaymath}
F_{n,t}(x) = \left( x-\sqrt{t} \right)^{n} + \left( x+\sqrt{t} \right)^{n}
\end{displaymath}
where $n \geq 4$ is a positive integer and $t$ is a negative integer.

One can find such examples for many different choices of
$\eta$ in Theorem~\ref{thm:gen-hypg1} below. Typically, we find that for fixed
$\eta$, there are infinitely many such examples with $n \leq 6$ and sometimes
some additional ones for larger $n$ too. The choice of $\eta$ here (essentially
$\sqrt{t}$) is unusual since for $n=6$, there are no such examples.

\section{Results}
\label{sect:results}

For positive integers $m$ and $n$ with $0<m<n$, $(m,n)=1$ and a non-negative
integer $r$, we put
\begin{displaymath}
X_{m,n,r}(x) = {} _{2}F_{1}(-r,-r-m/n;1-m/n;x),  
\end{displaymath}
where $_{2}F_{1}$ denotes the classical hypergeometric function. 

We let $D_{m,n,r}$ denote the smallest positive integer such that $D_{m,n,r} X_{m,n,r}(x) \in \bbZ[x]$.
For $d \in \bbZ$, we define $N_{d,m,n,r}$ to
be a positive integer such that $\left( D_{m,n,r}/ N_{d,m,n,r} \right)X_{m,n,r}\left( 1-\sqrt{d}\,z \right)
\in \bbZ \left[ \sqrt{\core(d)} \right] [z]$.
Here $\core(d)$ is the unique squarefree integer
such that $d/\core(d)$ is a square, with $\core(1)=1$. We will use
$v_{p}(x)$ to denote the largest power of a prime $p$ which divides
into the rational number $x$. We put
\begin{equation}
\label{eq:ndn-defn}
\cN_{d,n} =\prod_{p|n} p^{\min(v_{p}(d)/2, v_{p}(n)+1/(p-1))},
\end{equation}
and choose $\cC_{n}$ and $\cD_{n}$ such that
\begin{equation}
\label{eq:cndn-defn}
\max_{\stackrel{0<m<n/2}{\gcd(m,n)=1}} \left( \max \left( 1, \frac{\Gamma(1-m/n) \, r!}{\Gamma(r+1-m/n)},
\frac{n\Gamma(r+1+m/n)}{m\Gamma(m/n)r!} \right)
\frac{D_{m,n,r}}{N_{d,m,n,r}} \right)
< \cC_{n} \left( \frac{\cD_{n}}{\cN_{d,n}} \right)^{r}
\end{equation}
holds for all non-negative integers $r$, where $\Gamma(x)$ is the Gamma function.

One could choose $d, N_{d,m,n,r} \in \cO_{\bbK}$ such that
$\left( D_{m,n,r}/ N_{d,m,n,r} \right)X_{m,n,r}\left( 1-dx \right)
\in \cO_{\bbK} [x]$, with appropriate definitions of $\bbK$ and $\cN_{d,n}$.
However, our definition above avoids the required complications
and is sufficient for all our applications here.

\begin{theorem}
\label{thm:gen-hypg1}
Let $m$ and $n$ be as above, $t$, $u_{1}$ and $u_{2}$ be rational integers with
$t$ not a perfect square. Suppose that $\beta$ and $\gamma$ are algebraic integers
in $\bbQ \left( \sqrt{t} \right)$, with $\sigma$, the non-trivial element of
$\Gal \left( \bbQ \left( \sqrt{t} \right)/\bbQ \right)$. Put
\begin{eqnarray*}
\eta   & = & \left( u_{1} + u_{2} \sqrt{t} \right)/2, \\
\alpha & = & \frac{\beta(\eta/\sigma(\eta))^{m/n} \pm \sigma(\beta)}
                  {\gamma(\eta/\sigma(\eta))^{m/n} \pm \sigma(\gamma)},\\
g_{1}  & = & \gcd \left( u_{1}, u_{2} \right), \\
g_{2}  & = & \gcd \left( u_{1}/g_{1}, t \right),\\
g_{3}  & = & \left\{
			 \begin{array}{ll}
	             1 & \mbox{if $t \equiv 1 \bmod 4$ and $\left( u_{1}-u_{2} \right)/g_{1} \equiv 0 \bmod 2$}, \\
	             2 & \mbox{if $t \equiv 3 \bmod 4$ and $\left( u_{1}-u_{2} \right)/g_{1} \equiv 0 \bmod 2$},\\
	             4 & \mbox{otherwise,}
             \end{array}
             \right. \\
g      & = & g_{1}\sqrt{g_{2}/g_{3}}, \\
d      & = & \left( \eta-\sigma(\eta) \right)^{2}/g^{2} = u_{2}^{2}t/g^{2},
\end{eqnarray*}
\begin{eqnarray*}
E      & = & \frac{|g|\cN_{d,n}}{\cD_{n}\min \left( \left| u_{1} \pm \sqrt{u_{1}^{2}-u_{2}^{2}t} \right| \right)},\\
Q      & = & \frac{\cD_{n}\max \left( \left| u_{1} \pm \sqrt{u_{1}^{2}-u_{2}^{2}t} \right| \right)}{|g|\cN_{d,n}}, \\
\kappa & = & \frac{\log Q}{\log E} \mbox{ and } \\        
c      & = & 4 \sqrt{|2t|} \left( |\gamma| + |\sigma(\gamma)| \right) \cC_{n} Q \\
       &   & \times \left( \max \left( E, 5 \sqrt{|2t|} \left| 1- (\eta/\sigma(\eta))^{m/n} \right|
             |\beta-\alpha\gamma| \cC_{n}E \right) \right)^{\kappa},
\end{eqnarray*}
where the operation in the numerator of the definition of $\alpha$ matches the
operation in its denominator.

If $E > 1$ and either $($i$)$ $0 < \eta/\sigma(\eta) < 1$ or
$($ii$)$ $|\eta/\sigma(\eta)|=1$ with $\eta/\sigma(\eta) \neq -1$, then
$$
\left| \alpha - p/q \right| > \frac{1}{c |q|^{\kappa+1}} 
$$
for all rational integers $p$ and $q$ with $q \neq 0$. 
\end{theorem}

When $t$ is a perfect square, we have Corollary~2.6 in
\cite{Vout2}. Here too we can improve our choice of $d$ yielding the
following theorem.

\begin{theorem}
\label{thm:gen-hypg2}
Let $\bbK$ be an imaginary quadratic field and $m,n$ as above.
Let $a$ and $b$ be algebraic integers in $\bbK$ with the ideal $(a,b)=\cO_{\bbK}$
and either $a/b>1$ a rational number or $|a/b|=1$ with $a/b \neq -1$. Let
$\cC_{n}$, $\cD_{n}$ and $\cN_{d,n}$ be as above with $d=(a-b)^{2}$. Put
\begin{eqnarray*}
E      & = & \frac{\cN_{d,n}}{\cD_{n}} \left\{ \min \left( \left| \sqrt{a}-\sqrt{b} \right|, \left| \sqrt{a}+\sqrt{b} \right| \right) \right\}^{-2}, \\
Q      & = & \frac{\cD_{n}}{\cN_{d,n}} \left\{ \max \left( \left| \sqrt{a}-\sqrt{b} \right|, \left| \sqrt{a}+\sqrt{b} \right| \right) \right\}^{2}, \\
\kappa & = & \frac{\log Q}{\log E} \hspace{3.0mm} \mbox{ and }\\
c      & = & 4|a|\cC_{n}Q \left( 2.5\left| \frac{a(a-b)}{b} \right| \cC_{n}E \right)^{\kappa}.
\end{eqnarray*}

If $E > 1$, then 
$$
\left| (a/b)^{m/n} - p/q \right| > \frac{1}{c |q|^{\kappa+1}} 
$$
for all algebraic integers $p$ and $q$ in $\bbK$ with $q \neq 0$. 
\end{theorem}

In fact, in Theorem~3.2 and Theorem~3.5 of \cite{Vout3} we can take
$d=\left( \sigma(\eta)-\eta \right)/g$ and use the above definition of $\cN_{d,n}$.
In this way, the parameter $h$ that appears in both these theorems can also be
eliminated.

\subsection{New Irrationality Measures}

\begin{theorem}
\label{thm:n4gen}
Let $k=1$ or $3$. For a positive integer $b \geq 6$, write
$\left[ b\tan^{2}(k\pi/8) \right]=a_{1}a_{2}^{2}$, where $a_{1}$ is squarefree.
Suppose that $\gcd \left( a_{1}a_{2}^{2},b \right)=1$ and
$$
a_{1}a_{2}^{2}=b\tan^{2} (k\pi/8)+\epsilon,
$$
where $-0.5 < \epsilon < 0.5$. Let
\begin{eqnarray*}
\cN & = & \left\{
\begin{array}{cl}
1 & \mbox{if $a_{1}a_{2}b$ is even,} \\
4 & \mbox{if $a_{1}a_{2}b$ is odd and $a_{1} \equiv b \bmod 4$,} \\
8 & \mbox{if $a_{1}a_{2}b$ is odd and $a_{1} \not\equiv b \bmod 4$,}
\end{array}
\right. \\
\kappa & = & \left\{
\begin{array}{cl}
\displaystyle \frac{\log \left( 14.76b^{2}/\cN \right)}
{\log \left( \cN/(63.55\epsilon^{2}) \right)}
& \mbox{for $k=1$,} \vspace{1.0mm} \\
\displaystyle \frac{\log \left( 468.3b^{2}/\cN \right)}
{\log \left( \cN/(1.705\epsilon^{2}) \right)}
& \mbox{for $k=3$}
\end{array}
\right. \\
& & \text{ and } \\
c & = & (b/5) \left( 3 \cdot 10^{11} b^{3} \right)^{\kappa+1}.
\end{eqnarray*}

If the denominator of $\kappa$ is positive, then
\begin{equation}
\label{eq:result-n4gen}
\left| \sqrt{a_{1}b} \tan \left( \frac{k\pi}{8} \right) - \frac{p}{q} \right| > \frac{c}{|q|^{\kappa+1}} 
\end{equation}
for all integers $p$ and $q$ with $q \neq 0$.
\end{theorem}

\begin{note}
If $\epsilon=o(b^{-1/3})$, then this irrationality measure is better than the
Liouville bound. For example, the convergents, $a_{1}a_{2}^{2}/b$, in the
continued-fraction expansion of $\tan^{2}(\pi k/8)$ lead to such an improvement.

As in other applications of Thue's Fundamentaltheorem (e.g., \cite{Chen1, CV, LPV, TVW}),
where $\kappa$ approaches $1$ as a parameter like $b$ grows, here as $b$ in the
denominator of a continued-fraction convergent grows, $\kappa$ approaches $1$.
\end{note}

\begin{theorem}
\label{thm:n5gen}
Let $k=1$ or $2$. For a positive integer $b \geq 13$, write
$\left[ b\tan^{2}(k\pi/5) \right]=a_{1}a_{2}^{2}$, where $a_{1}$ is squarefree.
Suppose that $\gcd \left( a_{1}a_{2}^{2},b \right)=1$ and
$$
a_{1}a_{2}^{2}=b\tan^{2} (2k\pi/5)+\epsilon,
$$
where $-0.5 < \epsilon < 0.5$. With
$$
\cN_{1} = \left\{
\begin{array}{cl}
1 & \mbox{if $\gcd(5,a_{1}a_{2})=1$}, \\
5 & \mbox{if $5|a_{1}$}, \\
5^{5/4} & \mbox{if $5|a_{2}$},
\end{array}
\right.
$$
and

$$
\cN_{2} = \left\{
\begin{array}{cl}
1 & \mbox{if $a_{1}a_{2}b$ is even,} \\
4\sqrt{2} & \mbox{if $a_{1}a_{2}b$ is odd and $a_{1} \equiv b \bmod 4$,} \\
32 & \mbox{if $a_{1}a_{2}b$ is odd and $a_{1} \not\equiv b \bmod 4$,}
\end{array}
\right.
$$
let $\cN=\cN_{1}\cN_{2}$,
\begin{eqnarray*}
\kappa & = & \left\{
\begin{array}{cl}
\displaystyle \frac{\log \left( 5640b^{5/2}/\cN \right)}
{\log \left( \cN/(8.44b^{1/2} \epsilon^{2}) \right)}
& \mbox{for $k=1$,} \vspace{1.0mm} \\
\displaystyle \frac{\log \left( 48.26b^{5/2}/\cN \right)}
{\log \left( \cN/(57.68b^{1/2} \epsilon^{2}) \right)}
& \mbox{for $k=2$}
\end{array}
\right. \\
& & \text{ and } \\
c & = &(b/4000) \left( 8 \cdot 10^{14} b^{3} \right)^{\kappa+1}.
\end{eqnarray*}

If the denominator of $\kappa$ is positive, then
\begin{equation}
\label{eq:result-n5gen}
\left| \sqrt{a_{1}b} \tan \left( \frac{2k\pi}{5} \right) - \frac{p}{q} \right| > \frac{c}{|q|^{\kappa+1}} 
\end{equation}
for all integers $p$ and $q$ with $q \neq 0$.
\end{theorem}

\begin{note}
Here we require $\epsilon=o(b^{-2/3})$ to improve on the Liouville
irrationality measure. As above, all convergents, $a_{1}a_{2}^{2}/b$, in the
continued-fraction expansion of $\tan^{2}(2\pi k/5)$ lead to such an improvement.

However, unlike Theorem~\ref{thm:n4gen} and other applications of Thue's
Fundamentaltheorem, as $b$, in the denominator of a continued-fraction convergent,
grows, $\kappa$ approaches $5/3$.
\end{note}

\begin{theorem}
\label{thm:sporadic}
For all integers $p$ and $q$ with $q \neq 0$, we have
\noindent
\begin{equation}
\label{eq:n7t19}
\left| \sqrt{19} \tan \left( \frac{10\pi}{7} \right) - \frac{p}{q} \right| > 0.09 |q|^{-4.6},
\end{equation}
\begin{equation}
\label{eq:n7t39}
\left| \sqrt{39} \tan \left( \frac{8\pi}{7} \right) - \frac{p}{q} \right| > 0.007 |q|^{-3.28}, 
\end{equation}
\begin{equation}
\label{eq:n7t77}
\left| \sqrt{77} \tan \left( \frac{2\pi}{7} \right) - \frac{p}{q} \right| > 0.003 |q|^{-3.49}
\end{equation}
and
\begin{equation}
\label{eq:n13t7}
\left| \sqrt{7} \tan \left( \frac{18\pi}{13} \right) - \frac{p}{q} \right| > 0.02 |q|^{-5.68}.
\end{equation}
\end{theorem}

\section{Preliminary Results}

\subsection{Roots of $F_{n,t}(x)$}

The following lemma describes the roots.

\begin{lemma}
\label{lem:root-exp1}
Let $t$ be a negative integer.

\noindent
{\rm (i)} If $n$ is an odd positive integer, then the roots of $F_{n,t}(x)$
are $\sqrt{|t|} \tan (2k\pi/n)$ for $k=0,\ldots,n-1$.

\noindent
{\rm (ii)} If $n$ is an even positive integer, then the roots of $F_{n,t}(x)$
are $\sqrt{|t|} \tan ((2k+1)\pi/(2n))$ for $k=0,\ldots,n-1$.

\end{lemma}

\begin{proof}
Observe that
\begin{eqnarray*}
& & \left( \frac{\cos (\theta)}{\sqrt{|t|}} \right)^{n} F_{n,t}(\sqrt{|t|} \tan (\theta)) \\
& = & \left( \sin (\theta)-i\cos (\theta) \right)^{n} + \left( \sin (\theta)+i\cos (\theta) \right)^{n} \\
& = & \left( \cos (\theta-\pi/2) + i\sin (\theta-\pi/2) \right)^{n}
+ \left( \cos (\pi/2-\theta) + i\sin (\pi/2-\theta) \right)^{n} \\
& = & \cos (n(\theta-\pi/2)) + i\sin (n(\theta-\pi/2))
+ \cos (n(\pi/2-\theta)) + i\sin (n(\pi/2-\theta)) \\
& = & 2\cos (n(\pi/2-\theta)).
\end{eqnarray*}

(i) Letting $\theta=2k\pi/n$, we have
$$
n(\pi/2-\theta)=n(\pi/2-2k\pi/n)=n\pi/2-2k\pi.
$$

Since $n$ is odd, $2\cos(n\pi/2-2k\pi)=0$ and our result follows.

(ii) Here we let $\theta=(2k+1)\pi/(2n)$ and we find that
$$
n(\pi/2-\theta)=n(\pi/2-(2k+1)\pi/(2n))=n\pi/2-(2k+1)\pi/2.
$$

Since $n$ is even and $2k+1$ is odd, $2\cos(n\pi/2-(2k+1)\pi/2)=0$.
\end{proof}

\vspace{3.0mm}

The next lemma identifies roots of $F_{n,t}(x)$ with $\alpha$'s
in Theorem~\ref{thm:gen-hypg1}.

\begin{lemma}
\label{lem:alpha-exp}
Let $m=1$, $n$ as in Section~$\ref{sect:results}$, $t$ be a negative integer and
$z$ any integer. Put $\beta=\sqrt{t} \left( z+\sqrt{t} \right)$, $\gamma=z+\sqrt{t}$
and $\eta=\sqrt{t} \left(z-\sqrt{t} \right)^{n}$. Using subtraction in both the
numerator and denominator of the definition of $\alpha$ in Theorem~$\ref{thm:gen-hypg1}$,
we have
$$
\alpha =
\left\{
\begin{array}{ll}
\sqrt{|t|} \tan \left( \displaystyle \frac{\pi(n-\ell)}{2n} \right)      & \mbox{if $n$ is even} \\
\sqrt{|t|} \tan \left( \displaystyle \frac{2\pi((n-\ell)/4)}{n} \right)  & \mbox{if $n-\ell \equiv 0 \bmod 4$} \\
\sqrt{|t|} \tan \left( \displaystyle \frac{2\pi((3n-\ell)/4)}{n} \right) & \mbox{otherwise.}
\end{array}
\right.
$$
where $\left( z-\sqrt{t} \right)e^{\ell\pi i/n}/ \left( z+\sqrt{t} \right)$ is
the principal branch of $\left( \eta/\sigma \left( \eta \right) \right)^{1/n}$.
\end{lemma}

\begin{proof}
Substituting the values of $\beta$ and $\gamma$, we have
\begin{eqnarray*}
\alpha & = & \sqrt{t} \frac{\left( z + \sqrt{t} \right) \left( -\left( z-\sqrt{t} \right)^{n} /\left( z+\sqrt{t} \right)^{n} \right)^{1/n} + \left( z - \sqrt{t} \right)}
{\left( z + \sqrt{t} \right) \left( -\left( z-\sqrt{t} \right)^{n}/\left( z+\sqrt{t} \right)^{n} \right)^{1/n} - \left( z - \sqrt{t} \right)} \\
& = & \sqrt{t} \frac{(z-\sqrt{t})e^{\ell\pi i/n} + \left( z - \sqrt{t} \right)}
{(z-\sqrt{t})e^{\ell\pi i/n} - \left( z - \sqrt{t} \right)} \\
& = & \sqrt{t} \frac{e^{\ell\pi i/n} + 1}{e^{\ell\pi i/n} - 1}
= \sqrt{|t|} \frac{\sin(\ell\pi/n)}{1-\cos(\ell\pi/n)}
= \sqrt{|t|} \tan((n-\ell)\pi/(2n)),
\end{eqnarray*}
the last identity holds by a half-angle formula and symmetry about $\pi/2$.

Since we are taking an $n$-th root of $-1$ in $\left( \eta/\sigma \left( \eta \right) \right)^{1/n}$,
$\ell$ will be odd. If $n$ is even, then $n-\ell$ is odd and $\alpha$ is a root
of $F_{n,t}$.

If $n$ is odd, $n-\ell$ must be even. If $n-\ell \equiv 0 \bmod 4$, then
our result follows. Otherwise, notice that $\tan((n-\ell)\pi/(2n))
=\tan((3n-\ell)\pi/(2n))$ and $3n-\ell \equiv 0 \bmod 4$, completing our proof.
\end{proof}

\subsection{Arithmetic Estimates} 

\begin{lemma}
\label{lem:denom-est}
Let $\cC_{n}$ and $\cD_{n}$ be as defined in \eqref{eq:cndn-defn}.\\
{\rm (a)} For $n=4$, we can take $\cC_{n}=4.9 \cdot 10^{6}$ and $\cD_{n}=\exp(1.6)$.

\noindent
{\rm (b)} For $n=5$, we can take $\cC_{n}=8.8 \cdot 10^{9}$ and $\cD_{n}=\exp(1.37)$.

\noindent
{\rm (c)} For $n=7$, we can take $\cC_{n}=3.8 \cdot 10^{11}$ and $\cD_{n}=\exp(1.66)$.

\noindent
{\rm (d)} For $n=13$, we can take $\cC_{n}=1.9 \cdot 10^{13}$ and $\cD_{n}=\exp(2.21)$.
\end{lemma}

\begin{proof}
These are the $\cC_{1,n}$ values from Lemma~7.4(c) of \cite{Vout2} applied to these values of $n$.
\end{proof}

\begin{lemma}
\label{lem:g-and-n}
{\rm (a)} With $a_{1}, b$ and $\cN$ as in Theorem~$\ref{thm:n4gen}$, $g$ as in
Theorem~$\ref{thm:gen-hypg1}$ and $\cN_{d,4}$ as in \eqref{eq:ndn-defn},
$|g|\cN_{d,4}=2\cN a_{1}^{2}\sqrt{a_{1}b}$.

\noindent
{\rm (b)} With $a_{1}, b$ and $\cN$ as in Theorem~$\ref{thm:n5gen}$, $g$ as in
Theorem~$\ref{thm:gen-hypg1}$ and $\cN_{d,5}$ as in \eqref{eq:ndn-defn},
$|g|\cN_{d,5}=\cN a_{1}^{3}\sqrt{b}$.
\end{lemma}

\begin{proof}
(a) As we note in the proof of the Theorem~$\ref{thm:n4gen}$, we use
$z=a_{1}a_{2}$ and $t=-a_{1}b$, so
$u_{1} = 8a_{1}^{3}a_{2}b \left( a_{1}a_{2}^{2}-b \right)$ and
$u_{2} = 2a_{1}^{2}\left( a_{1}^{2}a_{2}^{4}-6a_{1}a_{2}^{2}b+b^{2} \right)$.

\noindent
$\bullet$ Determination of $g_{1}$

From the expressions for $u_{1}$ and $u_{2}$, we see that $2a_{1}^{2}|g_{1}$.
If $p>2$ is a prime dividing $g_{1}/\left( 2a_{1}^{2} \right)$, then either $p$
divides $a_{1}a_{2}b$ or else $a_{1}a_{2}^{2} \equiv b \bmod p$. The former
case is not possible since $a_{1}a_{2}$ and $b$ are relatively prime. In the
latter case, $p$ divides $4b^{2}$. But we have excluded $p=2$ and $p|b$ here.
Hence $g_{1}/\left( 2a_{1}^{2} \right)$ must be a power of two.

If one of $a_{1}a_{2}$ and $b$ is even and the other odd, then
$u_{2}/(2a_{1}^{2})=a_{1}^{2}a_{2}^{4}-6a_{1}a_{2}^{2}b+b^{2}$ is odd.
Hence $g_{1}/(2a_{1}^{2})$ is odd.

If $a_{1}a_{2}b$ is odd and $a_{1} \equiv b \bmod 4$, then 
$a_{1}^{2}a_{2}^{4}-6a_{1}a_{2}^{2}b+b^{2} \equiv 4 \bmod 8$. Therefore,
since $4|(u_{2}/(2a_{1}^{2}))$, $g_{1}/(8a_{1}^{2})$ is an odd integer.

If $a_{1}a_{2}b$ is odd and $a_{1} \not\equiv b \bmod 4$, then
$a_{1}^{2}a_{2}^{4}-6a_{1}a_{2}^{2}b+b^{2} \equiv 8 \bmod 16$.
Also $u_{2}/(-8a_{1}^{2}a_{2})=a_{1}a_{2}^{2}-b \equiv 2 \bmod 4$, so $g_{1}/(16a_{1}^{2})$
is an odd integer.

\noindent
$\bullet$ Determination of $g_{2}$

Since $2a_{1}^{2}|g_{1}$, we also have
$\gcd \left( u_{1}/g_{1}, t \right)
|\gcd \left( 4a_{1}a_{2}b \left( a_{1}a_{2}^{2}-b \right), a_{1}b \right)$.

Considering the cases examined for $g_{1}$, we find that $g_{2}=a_{1}b$.

\noindent
$\bullet$ Determination of $g_{3}$

Observe that
$$
\frac{u_{1}-u_{2}}{2a_{1}^{2}}
= -a_{1}^{2}a_{2}^{4}+6a_{1}a_{2}^{2}b-b^{2}+4a_{2}\left(a_{1}a_{2}^{2}-b \right) 
$$

If one of $a_{1}a_{2}$ and $b$ is even and the other is odd, then
$(u_{1}-u_{2})/g_{1}$ is odd and so $g_{3}=4$.

If $a_{1}a_{2}b$ is odd and $a_{1} \equiv b \bmod 4$, we saw above that
$\left( a_{1}^{2}a_{2}^{4}-6a_{1}a_{2}^{2}b+b^{2} \right)/4$ is odd. But
$a_{1}a_{2}^{2}-b$ is even.
Therefore $\left( u_{1}-u_{2} \right)/g_{1}$ is odd and $g_{3}=4$.

If $a_{1}a_{2}b$ is odd and $a_{1} \not\equiv b \bmod 4$, then $a_{1}^{2}a_{2}^{4}-6a_{1}a_{2}^{2}b+b^{2} \equiv 8 \bmod 16$.
Also $a_{1}a_{2}^{2}-b \equiv 2 \bmod 4$, so here $(u_{1}-u_{2})/g_{1}$ is even.
Furthermore, $t=-a_{1}b \equiv 1 \bmod 4$. Thus $g_{3}=1$.

\noindent
$\bullet$ Determination of $\cN_{d,4}$

We have
$$
d=\frac{u_{2}^{2}t}{g^{2}}
= \frac{g_{3}\left( a_{1}^{2}a_{2}^{4}-6a_{1}a_{2}^{2}b+b^{2} \right)^{2}}
{g_{1}^{2}/\left( 4a_{1}^{4} \right)}.
$$

To determine $\cN_{d,4}$ we need only consider the powers of $2$
dividing $d$.

If one of $a_{1}a_{2}$ and $b$ is even and the other is odd, then
$a_{1}^{2}a_{2}^{4}-6a_{1}a_{2}^{2}b+b^{2}$ and $g_{1}/(2a_{1}^{2})$ are odd
and $g_{3}=4$. Hence $2^{2} \parallel d$.

If $a_{1}a_{2}b$ is odd and $a_{1} \equiv b \bmod 4$, then
$\left( a_{1}^{2}a_{2}^{4}-6a_{1}a_{2}^{2}b+b^{2} \right)/4$ and $g_{1}/(8a_{1}^{2})$
are odd and $g_{3}=4$. Hence $2^{2} \parallel d$.

If $a_{1}a_{2}b$ is odd and $a_{1} \not\equiv b \bmod 4$, then
$a_{1}^{2}a_{2}^{4}-6a_{1}a_{2}^{2}b+b^{2} \equiv 8 \bmod 16$. Since $g_{1}/(16a_{1}^{2})$ is
also odd,
$\left( a_{1}^{2}a_{2}^{4}-6a_{1}a_{2}^{2}b+b^{2} \right) / \left( g_{1}/(2a_{1}^{2}) \right)$
is odd as well. Since $g_{3}=1$, we have $2^{0} \parallel d$.

Combining these observations, we have shown the following.

If one of $a_{1}a_{2}$ and $b$ is odd and the other is even, then
$|g|\cN_{d,4}=2a_{1}^{2}\sqrt{a_{1}b}$.

If $a_{1}a_{2}b$ is odd with $a_{1} \equiv b \bmod 4$, then
$|g|\cN_{d,4}=8a_{1}^{2}\sqrt{a_{1}b}$.

If $a_{1}a_{2}b$ is odd with $a_{1} \not\equiv b \bmod 4$,
then $|g|\cN_{d,4}=16a_{1}^{2}\sqrt{a_{1}b}$.

\vspace{3.0mm}

(b) The arguments to determine $g_{1}$, $g_{2}$ and $g_{3}$ are identical to those
for part~(a), so we only state the values of these quantities and
$\cN_{2}$ here.

If one of $a_{1}a_{2}$ and $b$ is odd and the other is even, then
$g_{1}=2a_{1}^{3}$, $g_{2}=b$ and $g_{3}=4$. So
$|g|=a_{1}^{3}\sqrt{b}$ and we can take $\cN_{2}=1$.

If $a_{1}a_{2}b$ is odd with $a_{1} \equiv b \bmod 4$, then
$g_{1}=8a_{1}^{3}$, $g_{2}=b$ and $g_{3}=2$. So
$|g|=4a_{1}^{3}\sqrt{2b}$ and we can take $\cN_{2} \geq 4\sqrt{2}$.

If $a_{1}a_{2}b$ is odd with $a_{1} \not\equiv b \bmod 4$, then
$g_{1}=32a_{1}^{3}$, $g_{2}=b$ and $g_{3}=1$. So $|g|=32a_{1}^{3}\sqrt{b}$
and we can take $\cN_{2} \geq 32$.

\noindent
$\bullet$ Determination of $\cN_{d,5}$

We have
$$
d=\frac{u_{1}^{2}t}{g^{2}}
= \frac{\sqrt{-g_{3}a_{1}}a_{2} \left( a_{1}^{2}a_{2}^{4} - 10a_{1}a_{2}^{2}b + 5b^{2} \right)}
{g_{1}/(2a_{1}^{2})}.
$$

To determine $\cN_{d,5}$ we are only interested the powers of $5$ dividing $d$.

If $5 \nmid a_{1}a_{2}$, then $5 \nmid d$ and $\cN_{d,5}=1$.

If $5|a_{2}$, then $25| \left( a_{2}\left( a_{1}^{2}a_{2}^{4} - 10a_{1}a_{2}^{2}b + 5b^{2} \right) \right)$,
and as we saw above $5 \nmid (g_{1}/(2a_{1}^{2}))$, so we can take $\cN_{d,5}=5^{5/4}$.

If $5|a_{1}$ and $5\nmid a_{2}$, then $5|| \left( a_{2}\left( a_{1}^{2}a_{2}^{4} - 10a_{1}a_{2}^{2}b + 5b^{2} \right) \right)$
and so $\cN_{d,5}=5$.

This argument justifies our choice of $\cN_{1}$ in Theorem~\ref{thm:n5gen}.
Combined with our results above about $\cN_{2}$, our lemma follows.
\end{proof}

\subsection{Analytic Estimates}

\begin{lemma}
\label{lem:anal-est}
(a) For any real $z$ with $-0.516<z<1$,
$$
1+z/2-z^{2}/8+z^{3}/16-z^{4}/16 \leq \sqrt{1+z}
\leq 1+z/2-z^{2}/8+z^{3}/16.
$$

(b) For any real $z$ with $0 \leq z \leq 0.62$,
$$
\arccos(1-z) \leq 1.5\sqrt{z}.
$$
\end{lemma}

\begin{proof}
(a) Using Maple, we find that
$$
\left( 1+z/2-z^{2}/8+z^{3}/16-z^{4}/16 \right)^{2}
- \left( 1+z \right)
= -\frac{3z^{4}}{64}-\frac{5z^{5}}{64}+\frac{5z^{6}}{256}
-\frac{z^{7}}{128}+\frac{z^{8}}{256}.
$$

The polynomial on the right-hand side has $z=0$, $z=-0.5161\ldots$ and $z=3$ as
its only real roots. This polynomial equals $-7/64$ at $z=1$ and $-7/65536$ at
$z=-1/2$. Therefore, it is at most zero for $-0.516<z<3$ and the desired lower bound
holds in this range.

A similar argument with
$$
\left( 1+z/2-z^{2}/8+z^{3}/16 \right)^{2}
- \left( 1+z \right)
= \frac{z^{4}}{64}-\frac{z^{5}}{64}+\frac{z^{6}}{256},
$$
shows that the polynomial on the right-hand side is non-negative for all real
$z$ and the desired upper bound holds in this range.

(b) $(d/dz) \arccos(1-z)=\left( 2z-z^{2} \right)^{-1/2}$, while
$(d/dz) 1.5\sqrt{z}=0.75z^{-1/2}$. For $0<z<1$, both of these derivatives are
positive and decreasing. The first one is less than the second one for $0<z<2/9$,
while the opposite is true for $2/9<z<1$. We find that
$\arccos(1-0.62)=1.1810\ldots$ and $1.5\sqrt{0.62}=1.1811\ldots$. Thus the
upper bound holds.
\end{proof}

\section{Proof of Theorems~\ref{thm:gen-hypg1} and \ref{thm:gen-hypg2}}

The arguments regarding $g_{1}, g_{2}$ and $g_{3}$ in Section~11 of \cite{Vout2}
continue to apply here. So Theorems~\ref{thm:gen-hypg1} and \ref{thm:gen-hypg2}
follow immediately from the following refinement of Lemma~7.4 of \cite{Vout2}.

\begin{lemma}
\label{lem:numer}
Suppose that $d$, $n$ and $r$ are non-negative integers with $d,n \geq 1$.
With $d_{1}=\gcd \left( d,n^{2} \right)$ and
$d_{2}=\gcd \left( d/d_{1}, n^{2} \right)$, we have
$$
\left( d_{1}^{\lfloor r/2 \rfloor} \prod_{p|d_{2}}p^{\min ( \lfloor rv_{p}(d_{2}) /2 \rfloor, v_{p}(r!))} \right) | N_{m,d,n,r}.
$$
\end{lemma}

\begin{proof}
This is a more general version of Proposition~5.1 of \cite{Chud} and we follow
the method of proof there. Using the reasoning there, we find that
\[
X_{m,n,r}\left( 1-\sqrt{d} \, x \right)
=\frac{r!}{(1-m/n)\cdots (r-m/n)} P_{r} \left( \sqrt{d} \, x \right),
\]
where
\[
P_{r}(x) = {2r \choose r} {} _{2}F_{1}(-r,-r-m/n;-2r;x).
\]
This is the $P_{-1}(x)$ in the proof of Proposition~5.1 of \cite{Chud}, but we
have added an index $r$.

So
$$
X_{m,n,r}\left( 1-\sqrt{d} \, x \right)
= \sum_{i=0}^{r} \left( \prod_{k=1}^{r-i} \frac{1}{kn-m} \right)
\frac{r! n^{r-i} d_{1}^{i/2}d_{2}^{i/2}d_{3}^{i/2}}{i!} { 2r-i \choose r} (-x)^{i},
$$
where $d_{3}=d/\left( d_{1}d_{2} \right)$. Since $(kn-m,n)=1$ for any integer
$k$, it is clear that $d_{1}^{\lfloor r/2 \rfloor}$ is a divisor of the
numerator of $X_{m,n,r}\left( 1-\sqrt{d} \, x \right)$.

Now suppose that $d_{2}>1$ and let $p$ be an odd prime divisor of $d_{2}$. Then
$p^{\lfloor i/2 \rfloor}/p^{v_{p}(i!)}$ is an integer, since $v_{p}(i!) \leq
i/(p-1) \leq i/2$. Hence we can remove a factor of $p^{v_{p}(r!)}$ from $r!$.
If $4|d_{2}$, then the same argument holds for $p=2$, while if $2 \parallel d_{2}$,
then we can remove a factor of $p^{\lfloor r/2 \rfloor}$. So in all cases, we
can remove a factor of
$p^{\min ( \lfloor rv_{p}(d_{2}) /2 \rfloor, v_{p}(r!))}$. Doing so for each
prime divisor of $d_{2}$ completes the proof.
\end{proof}

\section{Proof of Theorem~\ref{thm:n4gen}} 

We apply Theorem~\ref{thm:gen-hypg1} with $n=4$, $t=-a_{1}b$, $z=a_{1}a_{2}$,
$\beta=\sqrt{t} \left( z+\sqrt{t} \right)$, $\gamma=z+\sqrt{t}$ and
$\eta=\sqrt{t} \left( z-\sqrt{t} \right)^{n}$.

\subsection{Choice of $z$}

We check here that the above value of $z$ gives the algebraic numbers we
require. To do so, we find a sector containing $\left( z-\sqrt{t}\right) / \left( z+\sqrt{t}\right)$,
then use this to determine the principal branch of $\left( \eta /\sigma(\eta) \right)^{1/n}$
and hence $\ell$ in Lemma~\ref{lem:alpha-exp}.

We have
$$
\frac{z-\sqrt{t}}{z+\sqrt{t}}
= \frac{z^{2}+t-2z\sqrt{t}}{z^{2}-t}
= \frac{a_{1}a_{2}^{2}-b-2a_{2}\sqrt{-a_{1}b}}{a_{1}a_{2}^{2}+b}.
$$

We can write $a_{1}a_{2}^{2}-b=b\left( \tan^{2}(\pi k/8)-1 \right) +\epsilon$
and $a_{1}a_{2}^{2}+b=b\sec^{2}(\pi k/8)+\epsilon$, where $-0.5 < \epsilon < 0.5$.
So, with $b \geq 6$, for $k=1$, we have
$-0.838 < \Re \left( \left( z-\sqrt{t} \right)/ \left( z+\sqrt{t} \right) \right) < -0.593$.
Since $\Im \left( \left( z-\sqrt{t} \right) / \left( z+\sqrt{t} \right) \right)=-2a_{2} \sqrt{a_{1}b}/(a_{1}a_{2}^{2}+b)<0$,
\begin{equation}
\label{eq:n4-arg1}
-2.565 < \arg\left( \frac{z-\sqrt{t}}{z+\sqrt{t}} \right) < -2.2.
\end{equation}

Similarly, for $k=3$,
\begin{equation}
\label{eq:n4-arg2}
-0.8 < \arg \left( \frac{z-\sqrt{t}}{z+\sqrt{t}} \right) < -0.77.
\end{equation}

Next we bound the argument of $\left( \eta/\sigma(\eta) \right)^{1/4}$.

The real part of $\eta/\sigma(\eta)$ can be written as
$$
1 - \frac{2 \left( a_{1}^{2}a_{2}^{4}-6a_{1}a_{2}^{2}b+b^{2} \right)^{2}}
{ \left( a_{1}a_{2}^{2}+b \right)^{4} }
= 1- \frac{2\left( \left( a_{1}a_{2}^{2}-3b \right)^{2}-8b^{2} \right)^{2}}
{ \left( a_{1}a_{2}^{2}+b \right)^{4} },
$$
so we will show that this number, and hence $\eta/\sigma(\eta)$ itself, is near $1$.

Since $\tan^{4}(\pi k/8) - 6\tan^{2}(\pi k/8) +1=0$ and
$a_{1}a_{2}^{2}-3b = b \left( \tan^{2}(\pi k/8)-3 \right) + \epsilon$,
we have
$$
\left( a_{1}a_{2}^{2}-3b \right)^{2}-8b^{2}
= 2b\epsilon \left( \tan^{2}(\pi k/8)-3 \right) + \epsilon^{2}.
$$

So, for $k=1,3$ and $b \geq 6$,
$$
\left| 2b\epsilon \left( \tan^{2}(\pi k/8)-3 \right) + \epsilon^{2} \right|
< 5.75b|\epsilon|.
$$

Furthermore, for $b \geq 6$,
$$
1.088b < b\sec^{2}(\pi k/8)-0.5 < b\sec^{2}(\pi k/8)+\epsilon = a_{1}a_{2}^{2}+b.
$$

From the above expression for $\Re \left( \eta/\sigma(\eta) \right)-1$ and
these last two inequalities, we find that
$$
\left| \Re \left( \eta/\sigma(\eta) \right) - 1 \right|
< \frac{48\epsilon^{2}}{b^{2}}
$$
for $b \geq 6$. From Lemma~\ref{lem:anal-est}(b), we have
$$
|\arg\left( \eta/\sigma(\eta) \right)^{1/4})|<10.4|\epsilon|/(4b)<0.22.
$$

The interval $(-2.565+3\pi/4,-2.2+3\pi/4)$ is contained in the interval $(-0.22,0.22)$
while the interval $(-2.565+\pi/4,-2.2+\pi/4)$ does not intersect $(-0.22,0.22)$.
So from \eqref{eq:n4-arg1} and Lemma~\ref{lem:alpha-exp} with $\ell=3$, for
$k=1$, we find that $\alpha = \sqrt{|t|}\tan(\pi/8)$.

Considering \eqref{eq:n4-arg2} rather than \eqref{eq:n4-arg1},
$\alpha = \sqrt{|t|}\tan(3\pi/8)$ holds for $k=3$.

We also note here that from the above, for $b \geq 6$, we obtain
\begin{equation}
\label{eq:eta14-ub}
\left| \left( \eta/\sigma(\eta) \right)^{1/4} - 1 \right|
< \frac{2.6|\epsilon|}{b}.
\end{equation}

\subsection{Application of Theorem~\ref{thm:gen-hypg1}}

Since
$u_{1}^{2}-u_{2}^{2}t = 4 \left| \eta \right|^{2}
= 4a_{1}^{5}b \left( a_{1}a_{2}^{2}+b \right)^{4},
$
and $u_{1} = 8a_{1}^{3}a_{2}b \left( a_{1}a_{2}^{2}-b \right)$,
it follows that
$$
\frac{u_{1} \pm \sqrt{u_{1}^{2}-u_{2}^{2}t}}{2a_{1}^{2}\sqrt{a_{1}b}}
= 4\sqrt{a_{1}b} \, a_{2} \left( a_{1}a_{2}^{2}-b \right)
\pm \left( a_{1}a_{2}^{2}+b \right)^{2}.
$$

With $-0.5 < \epsilon < 0.5$, we have
\begin{eqnarray}
\label{eq:n4-v21}
\left( a_{1}a_{2}^{2}+b \right)^{2}
& = & b^{2}\sec^{4}(\pi k /8) + 2b\epsilon \sec^{2}(\pi k /8) + \epsilon^{2}, \\
a_{1}a_{2}^{2}b = b^{2}\tan^{2} (\pi k/8) + b\epsilon
& = & b^{2}\tan^{2} (\pi k/8) \left( 1 + \frac{\epsilon}{b\tan^{2} (\pi k/8)} \right), \nonumber \\
a_{1}a_{2}^{2} - b
& = & b \left( \tan^{2} (\pi k/8) - 1 \right) + \epsilon  \nonumber.
\end{eqnarray}

For $b \geq 6$ and $k=1$ or $3$, $|\epsilon/(b\tan^{2} (\pi k/8))|<0.49$, so
the bounds in Lemma~\ref{lem:anal-est}(a) apply and we have

\begin{eqnarray}
\label{eq:n4-v22LB}
& & \frac{\epsilon^{4}}{4b^{2}\tan^{7}(\pi k/8)}
- \frac{\epsilon^{5}}{4b^{3}\tan^{7}(\pi k/8)} \\
& < & 4\left( a_{1}a_{2}^{2}-b \right) \sqrt{a_{1}a_{2}^{2}b}
- \left\{ 4b^{2} \tan(\pi k/8) \left( \tan^{2} \left( \pi k/8 \right) - 1 \right) \right. \nonumber \\
& & \left. + 2b\epsilon \frac{3\tan^{2}(\pi k/8)-1}{\tan(\pi k/8)}
+ \frac{\epsilon^{2}}{2} \frac{3\tan^{2}(\pi k/8)+1}{\tan^{3}(\pi k/8)} \right. \nonumber \\
& & \left. - \frac{\epsilon^{3}}{4b}\frac{\tan^{2}(\pi k/8)+1}{\tan^{5}(\pi k/8)} \right\} \nonumber \\
\label{eq:n4-v22UB}
& < & \frac{\epsilon^{4}}{4b^{2}\tan^{5}(\pi k/8)}.
\end{eqnarray}

So, from \eqref{eq:n4-v21} and \eqref{eq:n4-v22LB}, and since
the left-hand side of \eqref{eq:n4-v22LB} is non-negative,
\begin{eqnarray}
\label{eq:n4-v2-+}
& & -4a_{2} \left( a_{1}a_{2}^{2}-b \right) \sqrt{a_{1}b}
+ \left( a_{1}a_{2}^{2}+b \right)^{2} \nonumber \\
& < & b^{2} \left( \sec^{4}(k\pi/8) - \left( 4\tan^{3}(k\pi/8) - 4\tan(k\pi/8) \right) \right) \nonumber \\
& + & 2b\epsilon \left( \sec^{2}(k\pi/8) - \frac{3\tan^{2}(k\pi/8)-1}{\tan(k\pi/8)} \right) \nonumber \\
& + & \frac{\epsilon^{2}}{2} \frac{2\tan^{3}(k\pi/8)-3\tan^{2}(k\pi/8)-1}{\tan^{3}(k\pi/8)}
+ \frac{\epsilon^{3}}{4b}\frac{\tan^{2}(k\pi/8)+1}{\tan^{5}(k\pi/8)}
\end{eqnarray}
and from \eqref{eq:n4-v21} and \eqref{eq:n4-v22UB},
\begin{eqnarray}
\label{eq:n4-v2++}
&   & 4a_{2} \left( a_{1}a_{2}^{2}-b \right) \sqrt{a_{1}b}
      + \left( a_{1}a_{2}^{2}+b \right)^{2} \nonumber \\
& < & b^{2} \left( \sec^{4}(k\pi/8) + \left( 4\tan^{3}(k\pi/8) - 4\tan(k\pi/8) \right) \right) \nonumber \\
&   & +2b\epsilon \left( \sec^{2}(k\pi/8) + \frac{3\tan^{2}(k\pi/8)-1}{\tan(k\pi/8)} \right) \nonumber \\
&   & +\frac{\epsilon^{2}}{2} \frac{2\tan^{3}(k\pi/8)+3\tan^{2}(k\pi/8)+1}{\tan^{3}(k\pi/8)}
- \frac{\epsilon^{3}}{4b}\frac{\tan^{2}(k\pi/8)+1}{\tan^{5}(k\pi/8)}.
\nonumber \\
&   & + \frac{\epsilon^{4}}{4b^{2}\tan^{5}(\pi k/8)}.
\end{eqnarray}

\subsubsection{$k=1$}

For $k=1$ and $b \geq 6$, $a_{1}a_{2}^{2}-b=b\left( \tan^{2}(\pi/8)-1 \right)+\epsilon=-0.8284\ldots b+\epsilon<0$.
Therefore, by substituting $k=1$ into \eqref{eq:n4-v2-+} and evaluating the trigonometric
functions, we obtain the upper bound
\begin{eqnarray}
\label{eq:n4k1MaxFcn}
& & \max \left| -4a_{2} \left( a_{1}a_{2}^{2}-b \right) \sqrt{a_{1}b}
\pm \left( a_{1}a_{2}^{2}+b \right)^{2} \right| \nonumber \\
& = &
-4a_{2} \left( a_{1}a_{2}^{2}-b \right) \sqrt{a_{1}b}
+ \left( a_{1}a_{2}^{2}+b \right)^{2} \nonumber \\
& < & b^{2} \left( 2.7451\ldots + \frac{4.6862\ldots\epsilon}{b}
- \frac{9.6568\ldots\epsilon^{2}}{b^{2}}
+ \frac{24.0208\ldots\epsilon^{3}}{b^{3}} \right).
\end{eqnarray}

If $\epsilon \leq 0$, then \eqref{eq:n4k1MaxFcn}
is at most $2.7451\ldots b^{2}$.
For $6 \leq b \leq 8$, $\epsilon<0$ and for $b=9$, $\epsilon=0.4558\ldots$,
so for $\epsilon>0$, we may assume $b \geq 9$. Now
$4.6862\ldots(\epsilon/b)
-9.6568\ldots \left( \epsilon/b \right)^{2}
+24.0208\ldots \left( \epsilon/b \right)^{3}<0.23465\ldots$ for
$\epsilon/b<0.5/9$ and hence the expression in \eqref{eq:n4k1MaxFcn} is at most
$2.9798b^{2}$. Thus
\begin{equation}
\label{eq:n4k1Max}
\max \left| -4a_{2} \left( a_{1}a_{2}^{2}-b \right) \sqrt{a_{1}b}
\pm \left( a_{1}a_{2}^{2}+b \right)^{2} \right|
< 2.9798b^{2},
\end{equation}

We turn now to the minimum. As above and applying \eqref{eq:n4-v2++}, we have
\begin{eqnarray}
\label{eq:n4k1MinFcn}
& & \min \left| -4a_{2} \left( a_{1}a_{2}^{2}-b \right) \sqrt{a_{1}b}
\pm \left( a_{1}a_{2}^{2}+b \right)^{2} \right| \nonumber \\
& < & \epsilon^{2} \left( 11.6568\ldots - \frac{24.0208\ldots\epsilon}{b}
+ \frac{20.503\ldots\epsilon^{2}}{b^{2}} \right).
\end{eqnarray}

If $\epsilon > 0$, then \eqref{eq:n4k1MinFcn} is at most
$11.6568\ldots \epsilon^{2}$.
As mentioned above, $\epsilon<0$ for $b=6$, $7, 8$ and $b \geq 12$. Calculating
\eqref{eq:n4k1MinFcn} directly for $b=6$, $7$
and $8$ and bounding it below by $\epsilon>-0.5$ for $b \geq 12$, we find that
\eqref{eq:n4k1MinFcn} is at most $12.83\epsilon^{2}$.
Hence, from Lemmas~\ref{lem:denom-est}(a) and \ref{lem:g-and-n}(a),
\begin{eqnarray*}
E & > & \frac{|g| \cN_{d,4}}{\cD_{4} 2a_{1}^{2}\sqrt{a_{1}b} \cdot 12.83\epsilon^{2}}
> \frac{\cN}{63.55\epsilon^{2}} \mbox{ and} \\
Q & < & \frac{\cD_{4}}{|g|\cN_{d,4}} 2a_{1}^{2}\sqrt{a_{1}b} \cdot 2.9798b^{2}
< \frac{14.76b^{2}}{\cN}.
\end{eqnarray*}

Finally, we determine an upper bound for $c$.
\begin{eqnarray*}
& & 4 \sqrt{|2t|} \left( |\gamma| + \left| \sigma(\gamma) \right| \right) \cC_{n} Q
\left( \max \left( E, 5 \sqrt{|2t|} \left| 1- \left( \eta/\sigma(\eta) \right)^{m/n} \right| |\beta-\alpha\gamma| \cC_{n}E \right) \right)^{\kappa} \\
& < & 8 \sqrt{2a_{1}b} \sqrt{a_{1}^{2}a_{2}^{2}+a_{1}b} \, 4.9\cdot 10^{6} \frac{14.76b^{2}}{\cN} \\
& & \times \left( 5\sqrt{2a_{1}b} \, \frac{2.6|\epsilon|}{b} \sqrt{a_{1}b}
\left| a_{1}a_{2}+\sqrt{-a_{1}b} \right| \left| 1 - i\tan \left( \frac{\pi}{8} \right) \right|
4.9\cdot 10^{6} \frac{\cN}{63.55\epsilon^{2}} \right)^{\kappa} \\
& < & \frac{8.2 \cdot 10^{8} a_{1}b^{5/2} \sqrt{a_{1}a_{2}^{2}+b}}{\cN}
\left( 1.54 \cdot 10^{6} a_{1}^{3/2}\sqrt{a_{1}a_{2}^{2}+b} \frac{\cN}{|\epsilon|} \right)^{\kappa} \\
& < & \frac{9.1 \cdot 10^{8} a_{1}b^{3}}{\cN}
\left( \frac{1.71 \cdot 10^{6} a_{1}^{3/2}b^{1/2}\cN}{|\epsilon|} \right)^{\kappa},
\end{eqnarray*}
since $a_{1}a_{2}^{2}+b=b\sec^{2}(\pi/8)+\epsilon<1.223b$ for $b \geq 6$ and
using \eqref{eq:eta14-ub}.

From $a_{1} \leq a_{1}a_{2}^{2}=b\tan^{2}(\pi/8)+\epsilon < 0.223b$ for $b \geq 6$,
we have
$$
c < \frac{2.1 \cdot 10^{8} b^{4}}{\cN}
\left( \frac{182,000 b^{2}\cN}{|\epsilon|} \right)^{\kappa}.
$$

The continued-fraction expansion of $\tan^{2}(\pi/8)$ is
$\left[ 0,5,\overline{1,4} \right]$. Using computation for small $q$ and the
fact that
\begin{equation}
\label{eq:cf-lb}
\frac{1}{\left( a_{i+1}+2 \right) q_{i}^{2}} < \left| \alpha-\frac{p_{i}}{q_{i}} \right|,
\end{equation}
where $a_{i+1}$ is the $i+1$-st partial fraction in the continued-fraction
expansion of $\alpha$ while $p_{i}/q_{i}$ is the $i$-th convergent, we find
that $|\epsilon|>1/(6b)$. Furthermore, since $\kappa>1$ and $\cN \leq 8$, have
$$
c < 3b \left( 9 \cdot 10^{6} b^{3} \right)^{\kappa+1}.
$$

\subsubsection{$k=3$}

Here we proceed in essentially the same way as for $k=1$, so we leave out
many of the details. By \eqref{eq:n4-v2-+} and \eqref{eq:n4-v2++}, we have
\begin{eqnarray}
\label{eq:n4k3Max}
&   & \max \left| -4a_{2} \left( a_{1}a_{2}^{2}-b \right) \sqrt{a_{1}b}
      \pm \left( a_{1}a_{2}^{2}+b \right)^{2} \right| \\
& = & 4a_{2} \left( a_{1}a_{2}^{2}-b \right) \sqrt{a_{1}b}
      + \left( a_{1}a_{2}^{2}+b \right)^{2} < 94.54b^{2} \mbox{ and} \nonumber \\
\label{eq:n4k3Min}
& & \min \left| -4a_{2} \left( a_{1}a_{2}^{2}-b \right) \sqrt{a_{1}b}
\pm \left( a_{1}a_{2}^{2}+b \right)^{2} \right|
< 0.3442\epsilon^{2}.
\end{eqnarray}

Hence, from Lemmas~\ref{lem:denom-est}(a) and \ref{lem:g-and-n}(a),

\begin{eqnarray*}
E & > & \frac{|g| \cN_{d,4}}{\cD_{4} 2a_{1}^{2} \sqrt{a_{1}b} \cdot 0.3442\epsilon^{2}}
> \frac{\cN}{1.705\epsilon^{2}}, \\
Q & < & \frac{\cD_{4}}{|g|\cN_{d,4}} 2a_{1}^{2} \sqrt{a_{1}b} \cdot 94.54b^{2}
< \frac{468.3b^{2}}{\cN} \mbox{ and} \\
c & < & (b/5) \left( 3 \cdot 10^{11} b^{3} \right)^{\kappa+1}.
\end{eqnarray*}

\section{Proof of Theorem~\ref{thm:n5gen}} 

We apply Theorem~\ref{thm:gen-hypg1} with $n=5$, $t=-a_{1}b$, $z=a_{1}a_{2}$,
$\beta=\sqrt{t} \left( z+\sqrt{t} \right)$, $\gamma=z+\sqrt{t}$ and
$\eta=\sqrt{t} \left( z-\sqrt{t} \right)^{n}$.

\subsection{Choice of $z$}

Again, the argument here is essentially the same as that used for the
choice of $z$ for Theorem~\ref{thm:n4gen}. For $b \geq 13$, we have
\begin{equation}
\label{eq:eta15-ub}
\left| \left( \eta/\sigma(\eta) \right)^{1/5} - 1 \right|
< \frac{1.1|\epsilon|}{b}.
\end{equation}

\subsection{Application of Theorem~\ref{thm:gen-hypg1}}

Here $u_{1}=2a_{1}^{3}b \left( 5a_{1}^{2}a_{2}^{4}-10a_{1}a_{2}^{2}b+b^{2} \right)$
and $u_{1}^{2}-u_{2}^{2}t = 4 \left| \eta \right|^{2} = 4a_{1}^{6}b \left( a_{1}a_{2}^{2}+b \right)^{5}$,
so
\begin{eqnarray*}
\frac{u_{1} \pm \sqrt{u_{1}^{2}-u_{2}^{2}t}}{2a_{1}^{3}\sqrt{b}}
& = & \left( 5 \left( a_{1}a_{2}^{2}-b \right)^{2} - 4b^{2} \right) \sqrt{b}
\pm \left( a_{1}a_{2}^{2}+b \right)^{2}\sqrt{a_{1}a_{2}^{2}+b}.
\end{eqnarray*}

We have
\begin{eqnarray}
\label{eq:n5-v2++}
& & \left( 5 \left( a_{1}a_{2}^{2}-b \right)^{2} - 4b^{2} \right) \sqrt{b}
+ \left( a_{1}a_{2}^{2}+b \right)^{2}\sqrt{a_{1}a_{2}^{2}+b} \nonumber \\
& < & 2b^{5/2}\sec^{5}(2\pi k/5)
+ 5b^{3/2}\epsilon \sec^{3}(2\pi k/5) + b^{1/2}\epsilon^{2} \left( 5 + \frac{15\sec(2\pi k/5)}{8} \right) \nonumber \\
& & + \frac{5\epsilon^{3}}{16b^{1/2}\sec(2\pi k/5)}
+ \frac{\epsilon^{5}}{16b^{5/2}\sec^{5}(2\pi k/5)}
\end{eqnarray}
and
\begin{eqnarray}
\label{eq:n5-v2-+}
& & \left( a_{1}a_{2}^{2}+b \right)^{2}\sqrt{a_{1}a_{2}^{2}+b}
- \left( 5 \left( a_{1}a_{2}^{2}-b \right)^{2} - 4b^{2} \right) \sqrt{b} \nonumber \\
& < & b^{1/2}\epsilon^{2} \left( 5-\frac{15\sec(2\pi k/5)}{8} \right)
-\frac{5\epsilon^{3}}{16b^{1/2}\sec(2\pi k/5)} \nonumber \\
& & +\frac{\epsilon^{4}}{16b^{3/2}\sec^{3}(2\pi k/5)}
+\frac{\epsilon^{5}}{16b^{5/2}\sec^{5}(2\pi k/5)}
+\frac{\epsilon^{6}}{16b^{7/2}\sec^{7}(2\pi k/5)}.
\end{eqnarray}

\subsubsection{$k=1$}
Applying the upper bound in \eqref{eq:n5-v2++}, for $b \geq 13$, we have
\begin{eqnarray*}
& & \max \left| \left( 5 \left( a_{1}a_{2}^{2}-b \right)^{2} - 4b^{2} \right) \sqrt{b}
\pm \left( a_{1}a_{2}^{2}+b \right)^{2}\sqrt{a_{1}a_{2}^{2}+b} \right| \\
& < & b^{5/2} \left( 709.77\ldots + \frac{169.44\ldots\epsilon}{b}
+ \frac{11.067\ldots\epsilon^{2}}{b^{2}}
+ \frac{0.096\ldots\epsilon^{3}}{b^{3}}
+ \frac{0.0001\ldots\epsilon^{5}}{b^{5}}
\right) \\
& < & 716.4b^{5/2}.
\end{eqnarray*}

Similarly, applying the upper bound in \eqref{eq:n5-v2-+},
\begin{eqnarray*}
& & \min \left| \left( 5 \left( a_{1}a_{2}^{2}-b \right)^{2} - 4b^{2} \right) \sqrt{b}
\pm \left( a_{1}a_{2}^{2}+b \right)^{2}\sqrt{a_{1}a_{2}^{2}+b} \right| \\
& < & b^{1/2}\epsilon^{2} \left( 1.0677\ldots + \frac{0.0966\ldots\epsilon}{b}
+ \frac{0.0002\ldots\epsilon^{3}}{b^{3}} \right) < 1.072b^{1/2}\epsilon^{2}.
\end{eqnarray*}

Hence, from Lemmas~\ref{lem:denom-est}(b) and \ref{lem:g-and-n}(b),
\begin{eqnarray*}
E & > & \frac{|g| \cN_{d,5}}{\cD_{5} 2a_{1}^{3}b^{1/2} \cdot 1.072b^{1/2}\epsilon^{2}}
> \frac{\cN}{8.44b^{1/2}\epsilon^{2}}, \\
Q & < & \frac{\cD_{5}}{|g|\cN_{d,5}} 2a_{1}^{3}b^{1/2} \cdot 716.4b^{5/2}
< \frac{5640b^{5/2}}{\cN} \mbox{ and } \\
c & < & (b/4000) \left( 8 \cdot 10^{14} b^{3} \right)^{\kappa+1}.
\end{eqnarray*}

\subsubsection{$k=2$}

As with $k=1$, for $b \geq 13$, we find that
\begin{eqnarray*}
& & \max \left| \left( 5 \left( a_{1}a_{2}^{2}-b \right)^{2} - 4b^{2} \right) \sqrt{b}
\pm \left( a_{1}a_{2}^{2}+b \right)^{2}\sqrt{a_{1}a_{2}^{2}+b} \right| \\
& < & b^{5/2} \left( 5.77\ldots + \frac{9.442\ldots\epsilon}{b}
+ \frac{2.682\ldots\epsilon^{2}}{b^{2}}
+ \frac{0.252\ldots\epsilon^{3}}{b^{3}}
+ \frac{0.021\ldots\epsilon^{5}}{b^{5}}
\right) \\
& < & 6.131b^{5/2}, \\
& & \min \left| \left( 5 \left( a_{1}a_{2}^{2}-b \right)^{2} - 4b^{2} \right) \sqrt{b}
\pm \left( a_{1}a_{2}^{2}+b \right)^{2}\sqrt{a_{1}a_{2}^{2}+b} \right| \\
& < & b^{1/2}\epsilon^{2} \left( 7.3176 \ldots + \frac{0.2528\ldots\epsilon}{b}
+ \frac{0.02166\ldots\epsilon^{3}}{b^{3}} \right)
< 7.328b^{1/2}\epsilon^{2}.
\end{eqnarray*}

Hence, from Lemmas~\ref{lem:denom-est}(b) and \ref{lem:g-and-n}(b),

\begin{eqnarray*}
E & > & \frac{|g| \cN_{d,5}}{\cD_{5} 2a_{1}^{3}b^{1/2} \cdot 7.328b^{1/2}\epsilon^{2}}
> \frac{\cN_{d,5}}{57.68b^{1/2}\epsilon^{2}}, \\
Q & < & \frac{\cD_{5}}{|g|\cN_{d,5}} 2a_{1}^{3}b^{1/2} \cdot 6.131b^{5/2}
< \frac{48.26b^{5/2}}{\cN} \mbox{ and } \\
c & < & (b/40000) \left( 8 \cdot 10^{14} b^{3} \right)^{\kappa+1}.
\end{eqnarray*}

\section{Larger $n$}

\subsection{Analysis} 

We can attempt the same proof for larger values of $n$.

For $n=6$, we just miss obtaining a theorem similar to Theorems~\ref{thm:n4gen}
and \ref{thm:n5gen}. For $k=1$ (the only $k$ we need consider for $n=6$),
\begin{eqnarray*}
\max \left( \left| \frac{u_{1} \pm \sqrt{u_{1}^{2}-u_{2}\sqrt{t}}}{g} \right| \right)
& < & b^{3} 2\sec^{6} \left( \frac{\pi}{12} \right), \\
\min \left( \left| \frac{u_{1} \pm \sqrt{u_{1}^{2}-u_{2}\sqrt{t}}}{g} \right| \right)
& < & 67.18 b \epsilon^{2}.
\end{eqnarray*}

Since $\tan^{2}(\pi/12)=1/(7+4\sqrt{3})$ is a quadratic irrational, $|\epsilon|>c_{1}/b$
for all positive integers, $b$. So even in the very best cases,
$$
\kappa= \frac{3\log(b)+c_{2}}{\log(b)+c_{3}},
$$
where $3c_{3}<c_{2}$ and hence $\kappa>3$.

Similarly, for larger values of $n$, we obtain
\begin{eqnarray*}
\max \left( \left| \frac{u_{1} \pm \sqrt{u_{1}^{2}-u_{2}\sqrt{t}}}{g} \right| \right)
& < & b^{n/2} c_{4}(n) \\
\min \left( \left| \frac{u_{1} \pm \sqrt{u_{1}^{2}-u_{2}\sqrt{t}}}{g} \right| \right)
& < & b^{n/2-2} \epsilon^{2}c_{5}(n).
\end{eqnarray*}

From Roth's theorem \cite{Roth}, $|\epsilon| < |b|^{-1-\delta}$ can only occur
finitely often for any $\delta>0$, so as $b$ grows, $\kappa$ approaches $n/(8-n)$.
Hence, for each $n \geq 7$, there are at most finitely many algebraic numbers
of the above form for which we can improve on Liouville's irrationality measure.

For $n \geq 9$, matters are even worse, since $n/2-2 > 2$, so, with at most
finitely many exceptions, we will not have $E>1$ and be unable to obtain any
irrationality measure from the hypergeometric method.

\subsection{Search Details} 

The algebraic numbers in Theorem~\ref{thm:sporadic} were found
by a computer search. The main idea behind the search is that $\eta/\sigma(\eta)$ must be near $1$
in order for us to be able to successfully apply the hypergeometric method.
This condition is the same as saying that $\eta-\sigma(\eta)=\sqrt{t}F_{n,t}(z)$
is small. That is, we choose $z$ near a root of $F_{n,t}$.

So for each $7 \leq n \leq 50$, our search was structured as follows.

(i) for each positive square-free integer $-1000 \leq t \leq -1$, and each integer $z$
from $\min_{F_{n,t}(\alpha)=0} \left( \sqrt{|t|} \alpha-10 \right)$ to
$\max_{F_{n,t}(\alpha)=0} \left( \sqrt{|t|} \alpha+10 \right)$, apply
Theorem~\ref{thm:gen-hypg1} to find values of $\kappa<\phi(n)-1$.

For smaller values of $t$, we observe that since $z$ is close to
$\sqrt{|t|} \tan(\theta)$ (for $\theta$ as in Lemma~\ref{lem:root-exp1}),
$z^{2}/|t|$ must be close to $\tan^{2}(\theta)$. As discussed in the previous
subsection, for larger $n$ we need the best approximations; and these
come from the continued-fraction expansion of $\tan^{2}(\theta)$. If $p/q$ is
a convergent in the continued-fraction expansion of $\tan^{2}(\theta)$ and we
write $p=p_{1} \cdot p_{2}^{2}$ where $p_{1}$ is a square-free integer, then we
can put $z=p_{1}p_{2}$ and $t=-p_{1}q$.

(ii) apply Theorem~\ref{thm:gen-hypg1} to $t$ and $z$ obtained
from the first 20 convergents in the continued-fraction expansion of the
appropriate $\tan^{2}(\theta)$'s.

The algebraic numbers in Theorem~\ref{thm:sporadic} were found
from step~(i). No further examples were found although there were
some near misses. The above calculations were
performed using PARI (version 2.3.3).

\section{Proof of Theorem~\ref{thm:sporadic}} 
\label{sec:8}

We will go through the details of the proof of \eqref{eq:n7t19},
identifying key quantities as we go along and then specifying the values of
these quantities for each of the remaining inequalities.

\subsection{Proof of \eqref{eq:n7t19}} 

We put $u_{1} = 2^{7} \cdot 13 \cdot 19^{4} \cdot 43$, $u_{2} = - 2^{7} \cdot 19^{4}$,
$m=1$, $n=7$, $t=-19$, $z=19$, $\beta=\sqrt{t} \left( z+\sqrt{t} \right)$ and
$\gamma=z+\sqrt{t}$. We have $\eta=\sqrt{t} \left( z-\sqrt{t} \right)^{n}$ and
$$
\frac{\eta}{\sigma(\eta)} = \frac{156231-559\sqrt{-19}}{156250}.
$$

Since we are using the principal branch when taking the $n$-th roots,
\begin{displaymath}
\left( \frac{156231-559\sqrt{-19}}{156250} \right)^{1/7}
= \frac{19-\sqrt{-19}}{19+\sqrt{-19}} e^{\pi i /7}.
\end{displaymath}

Thus we can apply Lemma~\ref{lem:alpha-exp} with $\ell=1$,
finding that $\alpha = \sqrt{19} \tan (10\pi /7)$.

\subsection{Application of Theorem~\ref{thm:gen-hypg1}}

Here $g_{1} = 2^{7} \cdot 19^{4}$ and $g_{2}= 1$.
Since $\left( u_{1}-u_{2} \right)/g_{1} \equiv 0 \bmod 2$ and $t \equiv 1 \bmod 4$,
we have $g_{3}=1$. Hence $g = 2^{7} \cdot 19^{4}$, $d=u_{2}^{2}t/g^{2}=-19$
and $\cN_{19,7}=1$. Also

\begin{eqnarray*}
\min \left( \left| u_{1} \pm \sqrt{u_{1}^{2}-u_{2}^{2}t} \right| \right)
& = & 2^{7} \cdot 19^{4} \left( -13 \cdot 43 + 2 \cdot 5^{3} \sqrt{5} \right) \mbox{ and} \\
\max \left( \left| u_{1} \pm \sqrt{u_{1}^{2}-u_{2}^{2}t} \right| \right)
& = & 2^{7} \cdot 19^{4} \left( 13 \cdot 43 + 2 \cdot 5^{3} \sqrt{5} \right).
\end{eqnarray*}

Thus, from Lemma~\ref{lem:denom-est}(c),
\begin{eqnarray*}
E & = & \frac{|g| \cN_{19,7}}{\cD_{7} \min \left( \left| u_{1} \pm \sqrt{u_{1}^{2}-u_{2}^{2}t} \right| \right)}
= 11.188347 \ldots,\\
Q & = & \frac{\cD_{7}}{|g| \cN_{19,7}} \max \left( \left| u_{1} \pm \sqrt{u_{1}^{2}-u_{2}^{2}t} \right| \right)
= 5879.998902 \ldots.
\end{eqnarray*}

So
$$
\kappa = \frac{\log Q}{\log E} < \frac{\log 5880}{\log 11.18834} < 3.59411,
$$
and
\begin{eqnarray*}
& & 4 \sqrt{38} \left( |\gamma| + |\sigma(\gamma)| \right) \cC_{7} Q \\
& & \times \left( \max \left( E, 5 \sqrt{38} \left| 1- (\eta/\sigma(\eta))^{1/7} \right|
             |\beta-\alpha\gamma| \cC_{7}E \right) \right)^{\kappa}
< 7 \cdot 10^{68}.
\end{eqnarray*}
Therefore, we can let $c=10^{69}$.

\subsection{Improved Constant}

By increasing $\kappa$ slightly, we can significantly reduce the size of $c$,
as in the proof of Corollary~2.2 of \cite{Vout1}.

We used Maple 8 to calculate the first $N=24,000$ partial fractions in the
continued-fraction expansion of $\sqrt{19} \tan (10\pi/7)$. This calculation
took 4750 seconds on a PC with an Intel Core i7-3630QM CPU running at 2.40 GHz.
The denominator of the $N=24,000$-th convergent is greater than $Q_{0}=10^{12000}$ 
and for all $q$ with $|q|>Q_{0}$,
\begin{displaymath}
\frac{10^{-69}}{|q|^{4.59411}} > \frac{0.09}{|q|^{4.6}}.
\end{displaymath}

The largest partial fraction found for $\sqrt{19} \tan (10\pi/7)$ was
$a_{1311}=21,976$. Applying this to \eqref{eq:cf-lb}, \eqref{eq:n7t19}
holds for $|q| \geq Q_{1}=19 > (0.09 \cdot (21976+2))^{(1/2.6)}$.
A direct check for all $|q| < Q_{1}$ completes the proof.

\subsection{Proof of \eqref{eq:n7t39}--\eqref{eq:n13t7}} 

As stated above, we proceed in the same way as for the proof of \eqref{eq:n7t19}
using the values in the accompanying table.

\begin{table}[ht]
\begin{center}
\begin{tabular}{||l|c|c|c||} \hline 
                     & \eqref{eq:n7t39} & \eqref{eq:n7t77} & \eqref{eq:n13t7} \\ \hline
$n$                  &    $7$                                 &  $7$    & $13$  \\ \hline
$t$                  &  $-39$                                 &  $-77$  & $-7$  \\ \hline
$z$                  &    $3$                                 &  $11$   &  $7$  \\ \hline
$u_{1}$              &  $2^{7}\cdot 3^{4}\cdot 13 \cdot 71$   &  $-2^{4}\cdot 7^{2} \cdot 11^{4} \cdot 167$ &
                        $-2^{13} \cdot 7^{7} \cdot 181$  \\ \hline
$u_{2}$              &  $-2^{7}\cdot 3^{4}$                   &  $2^{4} \cdot 11^{4}$  &
                        $-2^{13} \cdot 7^{7}$  \\ \hline
$\eta/\sigma(\eta)$  &  $\displaystyle\frac{32765 - 71 \sqrt{-39}}{32768}$ &
                        $\displaystyle\frac{4782958 - 1169\sqrt{-77}}{4782969}$ &
                        $\displaystyle\frac{16377 + 181\sqrt{-7}}{16384}$  \\ \hline
$\ell$               &  $5$ &  $3$  & $3$  \\ \hline
$g_{1}$              &  $2^{7} \cdot 3^{4}$   &  $2^{4} \cdot 11^{3}$  & $2^{13} \cdot 7^{7}$ \\ \hline
$g_{2}$              &  $13$   &  $7$   & $1$ \\ \hline
$g_{3}$              &  $1$    &  $2$   & $1$  \\ \hline
$d$                  &  $-3$   &  $-22$ & $-7$ \\ \hline
$\cN_{d,n}$          &  $1$   &  $1$    & $1$ \\ \hline
$E$                  &  $32.450014\ldots$   &  $75.606150 \ldots$  & $5.673393 \ldots$ \\ \hline
$Q$                  &  $2692.736355\ldots$   &  $46008.438040 \ldots$  & $3300.065595 \ldots$ \\ \hline
$\kappa$             &  $2.27$  & $2.4822$  &  $4.6675$ \\ \hline
$c$                  &  $7 \cdot 10^{48}$  & $2 \cdot 10^{54}$  &  $3 \cdot 10^{86}$ \\ \hline
$N$                  &  $10,000$  & $14,000$  &  $14,000$ \\ \hline
time(seconds)        &  $430$  & $980$  &  $1030$ \\ \hline
$Q_{0}$              &  $10^{5000}$  & $10^{7000}$  &  $10^{7000}$ \\ \hline
$\max a_{i}$         &  $a_{4021}=14,265$  & $a_{9118}=21,118$  &  $a_{2404}=303,427$ \\ \hline
$Q_{1}$              &  $37$  & $17$  &  $11$ \\ \hline
\end{tabular}                                                             
\end{center}
\end{table}

\begin{center}
{\sc Acknowledgements}
\end{center}

I thank my CM-friends and colleagues, Claude Levesque and Michel Waldschmidt,
for their encouragement to publish this work and the referee for their careful
reading of it.


\vspace{3.0mm}

\begin{thebibliography}{10}

\bibitem{Chen1}
Chen Jian Hua,   
A new solution of the Diophantine equation $X^{2}+1=2Y^{4}$, 
{\it J. Number Theory} {\bf 48} (1994), 62--74. 

\bibitem{CV}
Chen Jian Hua and P. M. Voutier,   
Complete solution of the diophantine equation $X^{2}+1=dY^{4}$ and a related family of quartic Thue equations, 
{\it J. Number Theory} {\bf 62} (1997), 71--99. 

\bibitem{Chud}
G. V. Chudnovsky,  
The method of Thue-Siegel,
{\it Annals of Math.} {\bf 117} (1983), 325--383. 

\bibitem{LPV}
G. Lettl, A. Peth\H{o} and P. M. Voutier,
Simple families of Thue inequalities,
{\it Trans. Amer. Math. Soc.} {\bf 351} (1999), 1871--1894.

\bibitem{Roth}
K. F. Roth,
Rational approximations to algebraic numbers, 
{\it Mathematika} {\bf 2} (1955), 1--20 and 168. 

%
\bibitem{TVW} 
A. Togb\'{e}, P. M. Voutier and P. G. Walsh,
Solving a family of Thue equations with an application to the equation $x^{2}-Dy^{4}=1$,
{\it Acta Arith.} {\bf 120} (2005), 39--58.

\bibitem{Vout1} 
P. M. Voutier, 
Rational approximations to $\sqrt[3]{2}$ and other algebraic numbers revisited, 
{\it Journal de Th\'{e}orie des Nombres de Bordeaux} {\bf 19} (2007), 263--288.

\bibitem{Vout2}
P. M. Voutier,   
Thue's Fundamentaltheorem, I: The General Case,
{\it Acta Arith.} {\bf 143} (2010), 101--144.

\bibitem{Vout3}
P. M. Voutier,   
Effective irrationality measures and approximations by algebraic conjugates,
{\it Acta Arith.} {\bf 149} (2011), 131--143.
\end{thebibliography}
\end{document}